\documentclass[12pt, a4Paper]{article}
\usepackage{cmap}     
\usepackage{mathtext}     
\usepackage[T2A]{fontenc}   
\usepackage[utf8]{inputenc}   
\usepackage[english]{babel} 
\usepackage{amsfonts}               
\usepackage{graphicx}               
\usepackage{amsmath}
\usepackage[nottoc]{tocbibind}
\usepackage{MnSymbol}%
\usepackage{wasysym}%
\graphicspath{{pictures/}}          
\DeclareGraphicsExtensions{.jpg}    
\title{\textbf{The Existence of Non-Equivariant Gromov Tori}}
\author{Grigoriy Yakovlev}
\date{}

\begin{document}

\maketitle

\begin{abstract}

In this paper, we address the following question: if a flat torus $\mathbb{T}^n$ is isometrically and minimally embedded into a sphere $\mathbb{S}^N$, must its translation group extend to the isometry group of the ambient sphere?

As shown by Robert Bryant, for $n=2$ the answer is positive. Furthermore, while Ying Lu, Peng Wang, and Zhenxiao Xie recently demonstrated that the answer is negative for immersions when $n \geq 3$, the question for embeddings remained open.

This problem is deeply tied to the work of Mikhail Gromov and Anton Petrunin concerning optimal curvature bounds. Petrunin proved that any immersion of a torus into a unit ball must have a maximum normal curvature of at least $\sqrt{\frac{3n}{n+2}}$. This bound is attained, for example, by families of tori constructed by Gromov. We call the tori that attain this optimal bound "Gromov tori".

In this work, we first demonstrate that any Gromov torus is intrinsically flat, lies within a sphere, and is minimal inside it. We then establish the necessary and sufficient conditions for defining these tori. Finally, we present our main result: for dimensions $n \ge 3$, there exists a non-equivariant embedded Gromov torus, which provides a definitive negative answer to the question above.

\end{abstract}

\tableofcontents
\newpage


\section{Introduction}
In this paper, we consider the following question:

\quad

\textbf{Question.} Let $\mathbb{T}^n$ be a flat torus isometrically embedded into a sphere $\mathbb{S}^N$ for some $N$. Suppose also that $\mathbb{T}$ is minimal in this sphere. Is it true that the translation group of $\mathbb{T}$ extends to the isometry group of the sphere?

\quad

A discussion on this topic can be found in \cite{MO}. As Robert Bryant showed in \cite{RB} for $n=2$, the answer is positive. He also conjectured that the answer is positive in the general case as well.

Also, in the case of immersions, it was shown in \cite{3.1} that for $n \geq 3$ the answer to the analogous question is negative, whereas the corresponding question for embeddings was explicitly left open.

The question itself was motivated by the works of Mikhail Gromov and Anton Petrunin. Mikhail Gromov in \cite{MG} constructed a family of geodesic equivariant subtori of the Clifford torus with normal curvatures equal to $\sqrt{\frac{3n}{n+2}}$ in all directions at each point. Anton Petrunin in \cite{TW} proved that this curvature estimate is optimal, namely, that any immersion of a torus into the unit ball must have a maximum normal curvature of at least $\sqrt{\frac{3n}{n+2}}$. The tori that achieve this bound are precisely what we call Gromov tori. Below we will show that any such torus carries a flat metric, lies in a sphere, and is minimal inside it.

In this work, we find the necessary and sufficient conditions for defining a Gromov torus, and we also show that for $n \geq 3$ there exists a non-equivariant embedded Gromov torus (i.e., one whose translation group does not extend to the isometry group of the sphere in the ambient space), which yields a negative answer to the posed question.

\section{Basic Concepts}
Let $L$ be a smooth $n$-dimensional manifold immersed into $\mathbb{R}^q$ with the induced Riemannian metric. Let $T_x$ and $N_x$ denote the tangent and normal spaces at the point $x \in L$, respectively.

The second fundamental form $\text{II}$ at the point $x$ is a symmetric quadratic form with arguments in $T_x$ and values in $N_x$. It is completely determined by the identity $\text{II}_x(v,v):=\gamma_v''(0)$, where $\gamma_v$ is a geodesic on $L$ such that $\gamma_v(0)=x$ and $\gamma_v'(0)=v$.

The normal curvature $k$ at the point $x$ in the direction $v$ is the value of the second form: $k_x(v):=\text{II}_x(v,v)$. We define the mean curvature vector $H(x) \in N_x$ as $H(x):=\sum_{i=1}^{n}{k_x(e_i)}$, where $\{e_i \}$ is some orthonormal basis in $T_x$.

Let $\text{Ж}(x)$ denote the average value of $|\text{II}_x|^2$ on the unit sphere in $T_x$, and let $\text{Sc}(x)$ denote the scalar curvature of $L$ at the point $x$.

A Gromov torus is a topological torus smoothly immersed into the unit ball $\mathbb{B}^q$ with a metric induced from $\mathbb{R}^{q+1}$, such that its normal curvatures are at most $\sqrt{\frac{3n}{n+2}}$ in all directions at every point.

\section{Properties of Gromov Tori}

We begin by stating some well-known results.

\textbf{Theorem 3.1 (Gauss Formula).} At any point $x \in L$, the following holds: 
$$\text{Sc}=\frac{3}{2}|H|^2-\frac{n\cdot(n+2)}{2}\cdot\text{Ж}$$

\quad

\textbf{Proof.} Consider an arbitrary point $p \in L$. First, assume that $\text{codim}(L)=1$. Let $k_1, \dots, k_n$ be the principal curvatures of $L$ at $p$. We know that $|H|^2=\sum_i{k_i^2}+2\sum_{i<j}{k_ik_j}$, and also that $n(n+2)\text{Ж}=3\sum_i{k_i^2}+2\sum_{i<j}{k_ik_j}$. The latter identity follows from the fact that $\text{Ж}$ is the average value of $(\sum_i{k_ix_i^2})^2$ on the unit sphere $\mathbb{S}^{n-1} \subset \mathbb{R}^n = T_p$, where $(x_1,\dots,x_n)$ are the standard coordinates in $\mathbb{R}^n$, and the functions $\frac{1}{3}n(n+2)x_i^4$ and $n(n+2)x_i^2x_j^2$ (for $i \neq j$) have a unit mean value.

By the standard Gauss formula, the scalar curvature is expressed through the principal curvatures as $\text{Sc}=2\sum_{i<j}{k_ik_j}$. Thus, $\text{Sc}$ is indeed expressed as $\frac{3}{2}|H|^2-\frac{n\cdot(n+2)}{2}\cdot\text{Ж}$.

If $\text{codim}(L)=k > 1$, then the second fundamental form at $p$ can be decomposed as a direct sum of real quadratic forms: $\text{II}_1 \oplus \dots \oplus \text{II}_k$, or $\text{II}=e_1\text{II}_1+\dots+e_k\text{II}_k$, where $\{ e_i \}$ is an orthonormal basis of $N_p$. From the case $k=1$, we have: $\text{Sc}_i=\frac{3}{2}|H_i|^2-\frac{n\cdot(n+2)}{2}\cdot\text{Ж}_i$ for $i=1,\dots,k$, where $\text{Sc}_i$, $H_i$, and $\text{Ж}_i$ are the quantities corresponding to $\text{II}_i$.

Since $H=\sum_i{H_i}$ and all $H_i$ are orthogonal to each other, $|H|^2=\sum_i{|H_i|^2}$. For similar reasons, $|\text{II}|^2=\sum_i{|\text{II}_i|^2}$, hence $\text{Ж}=\sum_i{\text{Ж}_i}$. Now we want to show that $\text{Sc}=\sum_i{\text{Sc}_i}$. Let $u$ and $v$ be orthonormal vectors in $T_p$. Then, by the standard Gauss formula, the sectional curvature in the direction $u \wedge v$ is expressed as $K_{u \wedge v}=\big< \text{II}(u,u),\text{II}(v,v) \big>-|\text{II}(u,v)|^2$, where $\text{Sc}$ is the sum of sectional curvatures over ordered, distinct $u$ and $v$ from the same orthonormal basis. Thus, it suffices to prove that $K_{u \wedge v}=\sum_i{K^i_{u \wedge v}}$. Let us expand $K_{u \wedge v}$:

$$K_{u \wedge v}=\big< \text{II}_1(u,u)e_1+\dots+\text{II}_k(u,u)e_k, \text{II}_1(v,v)e_1+\dots+\text{II}_k(v,v)e_k\big>-$$

$$-|\text{II}_1(u,v)e_1+\dots+\text{II}_k(u,v)e_k|^2=\sum_i{\big< \text{II}_i(u,u),\text{II}_i(v,v) \big>}-\sum_i{|\text{II}_i(u,v)|^2}=\sum_i{K^i_{u \wedge v}}$$

Finally, summing the expressions for all $\text{Sc}_i$, we obtain the required result.

\begin{flushright}
$\blacksquare$
\end{flushright}

\quad

\textbf{Corollary 3.2.} Let $L$ be closed. If $L$ is smoothly immersed in $\mathbb{B}^{q}$, then the average value of $|H|$ on $L$ is at least $n$.

If, in addition, $L$ is required to be smoothly immersed in $\mathbb{S}^q$, then $|H| \geq n$ at any point.

\quad

\textbf{Proof.} Consider the function $u:x \mapsto \frac{1}{2}|x|^2$ on $L$. Notice that $\Delta u(x)=n+\big< H(x),x \big>$. Then the average value of $\big< H(x),x \big>$ will be equal to $-n+ \frac{1}{|L|} \int_L{\Delta u(x)}dx=-n+\frac{1}{|L|}\int_L{\text{div}\nabla u(x)}dx=-n+\frac{1}{|L|}\int_{\partial L}{\nabla u(x)}dS(x)=-n$ by the Divergence Theorem (Ostrogradsky-Gauss formula), and since $\partial L=\varnothing$. Thus, if $|x| \leq 1$, the average value of $|H(x)|$ is indeed at least $n$.

If $L$ is immersed in a sphere, then $|x| \equiv 1$, and then the function $u$ under consideration becomes a constant. Then, at any point, $0=\Delta u(x)=n+\big< H(x),x \big>$, hence in this case we obtain $|H(x)| \geq n$ for any $x \in L$.

\begin{flushright}
$\blacksquare$
\end{flushright}

\quad

\textbf{Lemma 3.3.} Let $N$ be a smooth $m$-dimensional manifold immersed in $\mathbb{S}^{n-1} \subset \mathbb{R}^n$. Then if $|H_{N \subset \mathbb{R}^n}(p)|=m$ holds for $p \in N$, then $|H_{N \subset \mathbb{S}^{n-1}}(p)|=0$.

\quad

\textbf{Proof.} Let $\{ \gamma_i \}_{i=1}^m$ be a family of unit-speed geodesics on $N$ such that $\gamma_i(0)=p$, $\gamma'(0)=e_i$, where $\{ e_i \}$ is an orthonormal basis of $T_pN$; then $k_p(e_i)=\gamma''_i(0)$. Since $N \subset \mathbb{S}^{n-1}$, we have $|\gamma_i| \equiv 1$. Also, $\gamma''(0) \bot T_pN$ due to the unit-speed parametrization.

By assumption, we have $|\sum_{i=1}^m\gamma_i''(0)|=m$, and we want to prove that $\sum_{i=1}^m \text{pr}_{Tp\mathbb{S}^{n-1}} (\gamma_i''(0)) = \text{pr}_{Tp\mathbb{S}^{n-1}}(\sum_{i=1}^m\gamma_i''(0)) = 0$, where $\text{pr}_{Tp\mathbb{S}^{n-1}}(v)$ is the projection of $v$ onto $Tp\mathbb{S}^{n-1}$. 

By virtue of the unit-speed parametrization, we have $\big<\gamma_i'',\gamma_i\big>=\big<\gamma_i',\gamma_i\big>'-\big<\gamma_i',\gamma'_i\big>=-1$. Let $n:=\gamma_1(0)=\dots=\gamma_m(0)$ denote the position vectors of the geodesics at zero. Then we have $\big<\sum_{i=1}^m{\gamma_i++},n\big>=-m$, which implies $-m=|\sum_{i=1}^m{\gamma_i''}|\cos{\alpha}=m\cos{\alpha}$, where $\alpha$ is the angle between $n$ and $\sum_{i=1}^m{\gamma_i''}$. From this, we get $\cos{\alpha}=-1$, meaning that $n$ and $\sum_{i=1}^m{\gamma_i''}$ are parallel. But since $n$ is orthogonal to $T_p\mathbb{S}^{n-1}$ (as the curves $\gamma_i \subset \mathbb{S}^{n-1}$), $\sum_{i=1}^m{\gamma_i''}$ is also orthogonal to $T_p\mathbb{S}^{n-1}$. Thus, indeed, $\text{pr}_{Tp\mathbb{S}^{n-1}}(\sum_{i=1}^m\gamma_i''(0)) = 0$.

\begin{flushright}
$\blacksquare$
\end{flushright}

\quad

\textbf{Proposition 3.4.} Let $f:\mathbb{R}^n \rightarrow \mathbb{S}^q$ be a minimal isometric immersion. Then $\Delta f+nf = 0$.

\quad

\textbf{Proof.} Consider an arbitrary point $p \in f(\mathbb{R}^n)$. Then $p=f(x)=f(x_1,\dots,x_n)$ for some $x \in \mathbb{R}^n$.

Let $\{ x_i \}_{i=1}^n$ be the standard basis in $\mathbb{R}^n$. Since $f$ is an isometry, $\{ \frac{\partial f}{\partial x_i} \Big|_p \}_{i=1}^n$ is an orthonormal basis in $T_pf(\mathbb{R}^n)$. And since $f$ is a minimal immersion, $H_{f(\mathbb{R}^n) \subset \mathbb{S}^q}(p)=0$, meaning that the mean curvature vector $H_{f(\mathbb{R}^n) \subset \mathbb{R}^{q+1}}(p)$ must be orthogonal to $T_p\mathbb{S}^q$, and thus parallel to the position vector $p=f(x)$. We also know that $|f| \equiv 1$ since $f(\mathbb{R}^n) \subset \mathbb{S}^q$.

Now let us fix $x_2,\dots,x_n$ and consider $\gamma_1(t):=f(t,x_2,\dots,x_n)$ as a smooth curve with respect to the parameter $t$. Since at any point the partial derivatives form an orthonormal basis of the tangent space, $\gamma_1(t)$ is parameterized by arc length. Note that its curvature vector is then equal to $\gamma_1''(t)=\frac{\partial^2f}{\partial x_1^2} \Big|_{(t,x_2,\dots,x_n)}$. Since $\gamma_1 \subset \mathbb{S}^q$, we have $\big< \gamma_1,\gamma_1' \big> \equiv 0$, which yields $\big< \gamma_1,\gamma_1'' \big>=\big< \gamma_1,\gamma_1' \big>'-\big< \gamma_1',\gamma_1' \big>=-1$, i.e., $\big< f(t,x_2,\dots,x_n), \frac{\partial^2f}{\partial x_1^2} \Big|_{(t,x_2,\dots,x_n)} \big>=-1$.

Notice that $(t,x_2,\dots,x_n)$ is a geodesic in $\mathbb{R}^n$, as it is a parametrization of a coordinate line, and since $f$ is an isometry, $\gamma_1$ is a geodesic on $f(\mathbb{R}^n)$. Then, by the definition of normal curvature, we get that $k_{(t,x_2,\dots,x_n)}(\frac{\partial f}{\partial x_1} \Big|_{(t,x_2,\dots,x_n)})=\frac{\partial^2f}{\partial x_1^2} \Big|_{(t,x_2,\dots,x_n)}$.

Summing similar expressions for all $\gamma_i$ at the point $x=(x_1,\dots,x_n)$, we obtain $\big< f(x) , \Delta f(x) \big>=-n$. Note that $\Delta f(x) = \sum_{i=1}^n{\frac{\partial^2f}{\partial x_1^2} \Big|_{x}}=\sum_{i=1}^n{k_x(\frac{\partial f}{\partial x_1} \Big|_{x})}=H_{f(\mathbb{R}^n)\subset\mathbb{R}^{q+1}}(p) \parallel f(x)$. Thus, since $|f(x)| \equiv 1$, we finally find that $\Delta f(x)=-nf(x)$ at any point $x \in \mathbb{R}^n$.

\begin{flushright}
$\blacksquare$
\end{flushright}

\quad

In what follows, we will use the following propositions, the proofs of which can be found in \cite{TA} and \cite{MGBL}.

\quad

\textbf{Lemma 3.5.} Let $g$ be a metric on the torus $\mathbb{T}^n$, $n \geq 3$, let $\text{Sc}$ be its scalar curvature, and let $u:\mathbb{T}^n \rightarrow \mathbb{R}^n$ be a smooth positive function. Then $(\text{Sc} \cdot u - 4\cdot \frac{n-1}{n-2} \cdot \Delta u) \cdot u^{\frac{n-2}{n+2}}$ is the scalar curvature for the metric $u^{\frac{4}{n-2}} \cdot g$. As a consequence, the function $\text{Sc} \cdot u - 4\cdot \frac{n-1}{n-2} \cdot \Delta u$ takes a non-positive value at some point on the torus.

\quad

\textbf{Lemma 3.6.} Any Riemannian metric on $\mathbb{T}^n$ has a non-positive scalar curvature at some point. Moreover, if the scalar curvature on the torus is everywhere non-negative, then the torus is flat.

\quad

We will also use the proof of proposition $4.2$ from \cite{TW}.

\quad

\textbf{Proposition 3.7.} Any Gromov torus immersed in $\mathbb{B}^q$ lies in a sphere $\mathbb{S}^{q-1}$, is minimal in it, and has a flat metric.

\quad

\textbf{Proof.} Consider the function $u(x):=\exp{(-\frac{k}{2}|x|^2)}$ on our torus, where $k=\frac{3}{4} \cdot \frac{n-2}{n-1} \cdot n$. By Lemma 3.5, $\text{Sc}_1 := (\text{Sc} \cdot u - 4\cdot \frac{n-1}{n-2} \cdot \Delta u) \cdot u^{\frac{n-2}{n+2}}$ is the scalar curvature for some metric on $\mathbb{T}^n$. Let $\alpha(x):= \angle(H(x),x)$, $\beta(x):= \angle(N_x,x)$, and $r := |x|$. Then, according to the calculation carried out in \cite{TW}, $\text{Sc}_1$ can be rewritten as follows:

$$u^{\frac{n+2}{n-2}+1} \cdot \big( -\frac{n(n+2)}{2} \text{Ж} +\frac{3}{2} |H|^2 + 4\frac{n-1}{n-2}k|H|r\cos{\alpha}+4\frac{n-1}{n-2}(kn-k^2r^2\sin^2 {\beta}) \big)$$

Suppose there exists $x_0 \in \mathbb{T}^n$ such that $\text{Sc}_1(x_0) < 0$. As shown in \cite{TW}, the expression above is not less than $u^{\frac{n+2}{n-2}+1} \cdot (\frac{3}{2} n^2-\frac{n(n+2)}{2}\text{Ж})$. It follows that $\frac{n(n+1)}{2}\cdot \text{Ж}(x_0) > \frac{3n}{n+2}$, i.e., $\text{Ж}(x_0) > \frac{3n}{n+2}$, which contradicts the definition of a Gromov torus. Then, by Lemma 3.6, we have that $\mathbb{T}^n$ is a flat torus with respect to the metric $u^{\frac{4}{n-2}} \cdot g$, meaning $\text{Sc}_1 \equiv 0$. This implies that $\text{Ж} \geq \frac{3n}{n+2}$ at all points of the torus, hence $\text{Ж} \equiv \frac{3n}{n+2}$. Consequently, all normal curvatures are equal to $\sqrt{\frac{3n}{n+2}}$.

Notice that $\cos^2{\alpha} + \sin^2{\beta} \leq 1$. At the same time, $1 < \frac{3}{2}\frac{n-2}{n-1} < 2$, $r^2 \leq \cos^2{\alpha}$ in the case $n \geq 5$ (by proposition $4.1$ from \cite{TW}) and $\frac{3}{2} \geq \frac{9}{4}\frac{n-1}{n-2}$ in the case $n \leq 4$. Let us consider two cases:

1) $n \leq 4$. Then $\frac{3}{2} |H|^2 + 4\frac{n-1}{n-2}k|H|r\cos{\alpha}+4\frac{n-1}{n-2}(kn-k^2r^2\sin^2 {\beta}) \geq \frac{3}{2} n^2+\frac{3}{2} n^2(1-r^2)+\frac{3}{2}(|H|+nr\cos{\alpha})^2 \geq \frac{3}{2} n^2$. From the arguments above, we see that these inequalities must turn into equalities at all points. In particular, this implies that $r \equiv 1$ and $|H|+nr\cos{\alpha} \equiv 0$. Consequently, $\mathbb{T}^n \subset \mathbb{S}^{q-1}$, as well as $|H| \leq n$. But by Corollary 3.2, we know that then $|H| \geq n$ at all points, meaning that in fact $|H| \equiv n$.

2) $n \geq 5$. Then $\frac{3}{2} |H|^2 + 4\frac{n-1}{n-2}k|H|r\cos{\alpha}+4\frac{n-1}{n-2}(kn-k^2r^2\sin^2 {\beta}) \geq \frac{3}{2} n^2+\frac{3}{2} n^2(\sin^2{\beta}\cdot(2-\frac{3}{2}\frac{n-2}{n-1})+\sin^4{\beta} \cdot(\frac{3}{2}\frac{n-2}{n-1}-1))+\frac{3}{2}(|H|+nr\cos{\alpha})^2 \geq \frac{3}{2}n^2$. Analogously to the first case, all inequalities must turn into equalities. This is possible only when $|H|+nr\cos{\alpha}\equiv 0$ and $\sin{\beta}\equiv 0$. By Corollary 3.2, the average value of $|H| \geq n$, which implies $|H| \equiv n$, from which it follows that $|x| \equiv 1$. This means that $\mathbb{T}^n \subset \mathbb{S}^{q-1}$.

Thus, in both cases, we obtained that $\mathbb{T}^n \subset \mathbb{S}^{q-1}$ and $|H| \equiv  n$. Then, by Lemma 3.3, this means that $\mathbb{T}^n$ is minimal in the sphere. In turn, it follows from the Gauss formula that $|\text{Sc}| \equiv 0$. Hence, by Lemma 3.6, $\mathbb{T}^n$ is flat.

\begin{flushright}
$\blacksquare$
\end{flushright}

\quad

\section{Fourier Analysis on the Torus}

Consider the space $\mathbb{R}^n$ with the standard metric, and a lattice $A \mathbb{Z}^n$, where $A$ is some non-singular matrix of size $n \times n$. The flat torus $\mathbb{T}^n$ as a topological space coincides with $\mathbb{R}^n / A \mathbb{Z}^n$, and the metric is induced via factorization from $\mathbb{R}^n$.

Let us define a family of smooth functions on the torus by $l_k(x) := \big<x,A^{-T}k \big>$, where $k \in \mathbb{Z}^n$. We also define a family of smooth maps $u_k: \mathbb{T}^n \rightarrow \mathbb{S}^1 \subset \mathbb{C}$ by the formula $u_k := e^{2\pi i l_k}$.

This section outlines well-known facts, an alternative proof of which can be found in \cite{FT}.

\quad

\textbf{Theorem 4.1.} $\{ u_k \}_{k \in \mathbb{Z}^n}$ form a Fourier basis for $L^2(\mathbb{T}^n)$.

\quad

\textbf{Proof.} First, let us prove that $\{ u_k \}$ is an orthogonal system. To do this, we transform the following expression:

$$\int_{\mathbb{T}^n}{e^{2\pi i l_k}}dx=\int_{\mathbb{T}^n}{e^{2\pi i \big<x,A^{-T}k \big>}}dx=\int_{\mathbb{T}^n}{e^{2\pi i \big<A^{-1}x,k \big>}}dx=|A|\int_{A^{-1}\mathbb{T}^n}{e^{2\pi i \big<x,k \big>}}dx=|A|\int_{[0,1]^n}{e^{2\pi i \big<x,k \big>}}dx$$

Now, let $k_1, k_2\in\mathbb{Z}^n$. Then:

$$\frac{1}{|A|}\int_{\mathbb{T}^n}{e^{2\pi i l_{k_1}} \overline{e^{2\pi i l_{k_2}}}}dx=\frac{1}{|A|}\int_{\mathbb{T}^n}{e^{2\pi il_{k_1-k_2}} }dx= \int_{[0,1]^n}{e^{2\pi \big< x,k_1-k_2 \big>} }dx =$$
 
 $$=\int_{[0,1]}{e^{2\pi i (k_1^n-k_2^n)x_n}}\dots\int_{[0,1]}{e^{2\pi i (k_1^1-k_2^1)x_1}}dx_1\dots dx_n$$

If $k \neq l$, we can assume that they do not coincide in the first coordinate. Then:

$$\int_{\mathbb{T}^n}{e^{2\pi i \big< x,k \big>} \overline{e^{2\pi i \big< x,l \big>}}}dx=\int_{[0,1]}{e^{2\pi i (k_1^n-l_2^n)x_n}}\dots\bigg(\frac{e^{2\pi i (k_1^1-k_2^1)x_1}}{2\pi i (k_1^1-k_2^1)} \bigg|_{0}^1\bigg)\dots dx_n=0$$

The last step holds because $k_1^1-k_2^1 \in \mathbb{Z}$. This means $\big< u_{k_1},u_{k_2} \big>_{L_2}=0$ when $k_1 \neq k_2$, so it is indeed an orthogonal system.

Now we show that this system actually forms a basis. To this end, we wish to apply the Stone-Weierstrass theorem for the algebra $\text{Span} (\{ u_k \}_{k \in \mathbb{Z}^n})=:U$. We check the conditions:

\begin{itemize}
\item $1=e^{2\pi i \big< x,0 \big>}=u_0 \in U$.
\item Let $x, y$ be distinct points on the torus; we wish to separate them with a function from $U$. Consider $u_k(x) u_k(y)^{-1} = e^{2\pi i \big< x-y, A^{-T}k \big>}$. If this expression were equal to $1$ for all $k \in \mathbb{Z}^n$, it would follow that $\big< x-y, A^{-T}k \big> \in \mathbb{Z}$ for all $k \in \mathbb{Z}^n$. This is equivalently written as $\big< A^{-1}(x-y), k \big> \in \mathbb{Z}$ for all $k \in \mathbb{Z}^n$, which implies $A^{-1}(x-y) \in \mathbb{Z}^n$. Hence, $x=y$ on the torus, yielding a contradiction.

\item If $f \in U$, then $f=\sum_{j=1}^m{a_k u_k}$, and then $\overline{f}=\sum_{j=1}^m{\overline{a_k} \overline{u_k}}=\sum_{j=1}^m{\overline{a_k} u_{-k}} \in U$.
\end{itemize}

Thus, we can indeed apply the Stone-Weierstrass theorem and conclude that $U$ is dense in $C(\mathbb{T}^n)$. Since $C(\mathbb{T}^n)$ is dense in $L^2(\mathbb{T}^n)$, it follows that $U$ is dense in $L^2(\mathbb{T}^n)$ as well. Therefore, no non-zero element in $L^2(\mathbb{T}^n)$ is orthogonal to the entire algebra $U$, which implies that functions $\{ u_k \}_{k \in \mathbb{Z}^n}$ form a basis for $L^2(\mathbb{T}^n)$.

\begin{flushright}
$\blacksquare$
\end{flushright}

\quad

Now let $f: \mathbb{T}^n \rightarrow \mathbb{C}$ be a smooth function. It follows from the theorem that it can be expanded into a Fourier series: $f(x)=\sum_{k \in \mathbb{Z}^n}{a_k\cdot u_k(x)}$.

\quad

\textbf{Corollary 4.2.} If $f$ is such that $\Delta f+nf = 0$, then there is a finite number of non-zero coefficients in its Fourier series expansion.

\quad

\textbf{Proof.} Let $A^{-T}k\in A^{-T}\mathbb{Z}^n$ be expanded in the standard basis $\{ \varepsilon_j \}$ as $\sum_{l=1}^n{\varepsilon_lk_l}$. Consider the Laplace operator applied to a Fourier basis element $u_k$:

$$\Delta u_k=\sum_{j=1}^n{\frac{\partial^2e^{2\pi i \big< x,A^{-T} k \big>}}{\partial x_j^2}}=-4\pi^2\bigg(\sum_{j=1}^n{k_j^2}\bigg)e^{2\pi i \big< x,A^{-T}k \big>}=-4\pi^2|A^{-T}k|^2e^{2\pi i \big< x,A^{-T}k \big>}=-4\pi^2|A^{-T}k|^2u_k$$

Since $f(x)=\sum_{k \in \mathbb{Z}^n}{a_k\cdot u_k(x)}$, from the equation $\Delta f+nf = 0$, for each $k \in \mathbb{Z}^n$ such that $a_k \neq 0$, we get: $4\pi^2|A^{-T}k|^2=n$, which means $|A^{-T}k|=\frac{\sqrt{n}}{2\pi}$. The left-hand side of this expression is the length of some vector from the lattice $A^{-T}\mathbb{Z}^n$. Obviously, there is a finite number of vectors of length $\frac{\sqrt{n}}{2\pi}$ in such a lattice, which implies that the number of non-zero $a_k$ must also be finite.

\begin{flushright}
$\blacksquare$
\end{flushright}

\quad

\section{Non-Equivariant Gromov Torus}

From the preceding arguments, it follows that any Gromov torus is defined by a smooth map $f: \mathbb{R}^{n} \rightarrow \mathbb{C}^{q+1}$ with the following properties: $\text{Im} f \subset \mathbb{R}^{q+1}$, $\Delta f + nf=0$, $|f| \equiv 1$, $|D_uf| \equiv 1$, $|D_u^2f| \equiv \sqrt{\frac{3n}{n+2}}$ (where $u$ is an arbitrary unit vector in $\mathbb{R}^n$), and there exists a lattice $\Lambda \subset \mathbb{R}^n$ such that $f(x+a)=f(x)$ for any $a \in \Lambda$.

We also established that any such map can be represented as a finite sum:

$$f(x)=\sum_{k}{c_k \cdot \exp{(2\pi i \big< k,x \big>)}}$$

Where the sum is taken over all $k$ from the dual lattice $\Lambda^{*}$ (the lattice corresponding to the inverse transpose matrix of the original lattice) satisfying the condition $|k| = \frac{\sqrt{n}}{2\pi}$, and $\{ c_k \}$ is a set of vectors in $\mathbb{C}^{q+1}$.

In fact, the converse is also true:

\quad

\textbf{Lemma 5.1.} Let $f: \mathbb{R}^{n} \rightarrow \mathbb{C}^{q+1}$ be a smooth map and let $\Lambda \subset \mathbb{R}^n$ be a lattice satisfying the following conditions:

\begin{enumerate}
\item $f(x)=\sum_{k}{c_k \cdot \exp{(2\pi i \big< k,x \big>)}}$, where the sum is taken over all $k$ from the lattice $\Lambda^{*}$ satisfying $|k| = \frac{\sqrt{n}}{2\pi}$, and $\{ c_k \}$ is a set of vectors in $\mathbb{C}^{q+1}$;
\item $\text{Im} f \subset \mathbb{R}^{q+1}$;
\item $|f| \equiv 1$;
\item $|D_{e_i}f| \equiv 1$, where $\{{e_i}\}$ is the standard basis in $\mathbb{R}^n$;
\item $\big<D_{e_i}f,D_{e_j}f\big> \equiv 0$ for $i \neq j$;
\item $|D_{e_i}^2f| \equiv \sqrt{\frac{3n}{n+2}}$;
\item $|D^2_{\frac{e_i+e_j}{\sqrt{2}}}f|=\sqrt{\frac{3n}{n+2}}$ for $i \neq j$.
\end{enumerate}

Then $f$ defines a Gromov torus with lattice $\Lambda$.

\quad

\textbf{Proof.} Note that for any $x \in \mathbb{R}^n$ and $a \in \Lambda$, we have $f(x+a)=\sum_{k}{c_k \cdot \exp{(2\pi i \big< k,x+a \big>)}}=\sum_{k}{c_k \cdot \exp{(2\pi i \big< k,x \big>)}\cdot\exp{(2\pi \big< k,a \big>)}} = f(x)$, because $\big<k,a \big>$ is an integer. From conditions 2 and 3, it follows that $\text{Im}f \subset \mathbb{S}^q$. Then from conditions 4 and 5, it follows that $f$ is an isometric immersion of $\mathbb{R}^n/\Lambda \mathbb{Z}^n$ into $\mathbb{S}^{q}$. Thus, the image of any line will be a geodesic in the immersed torus, and hence, by property 6, all normal curvatures of the torus in the directions $df(e_i)$ will be equal to $\sqrt{\frac{3n}{n+2}}$. It is easy to see that by property 7, $\text{II}(dfe_i,dfe_j)=0$ for $i \neq j$, which implies that all normal curvatures will be equal to $\sqrt{\frac{3n}{n+2}}$. Thus, we find that $f$ defines a Gromov torus with lattice $\Lambda$.

\begin{flushright}
$\blacksquare$
\end{flushright}

\quad

\textbf{Lemma 5.2.} Let $f(x)=\sum_{k}{c_k \cdot \exp{(2\pi i \big< k,x \big>)}}$ define an embedded Gromov torus. It is claimed that this torus is equivariant (i.e., the translation group of this torus extends to the isometry group of the sphere in which it is embedded) if and only if $c_k \perp c_l$ for $k \neq l$.

\quad

\textbf{Proof.} We will prove the lemma according to the following scheme:

1) Torus is equivariant $\Rightarrow$ torus is an orbit of a torus in $O(q)$.

2) Torus is an orbit of a torus in $O(q)$ $\Rightarrow$ Torus is equivariant.

3) Torus is an orbit of a torus in $O(q)$ $\Rightarrow$ $c_k \perp c_j$ for $k \neq j$.

4) $c_k \perp c_j$ for $k \neq j$ $\Rightarrow$ torus is an orbit of a torus in $O(q)$.

\quad

Step 1. Let $G$ be the translation group of the torus $\mathbb{T}^n$, then $G = \mathbb{T}^n$, and the action is defined as $a(f(x)):=f(a+x)$, $a,x \in \mathbb{T}^n$. By hypothesis, we have an injective homomorphism $\varphi:\mathbb{T}^n \rightarrow O(q)$ such that $f(x+a)=\varphi(a)f(x)$. This means $\text{Im}f=\text{Orb}_{\varphi(\mathbb{T}^n)}(f(0))$. Moreover, $\varphi(\mathbb{T}^n)$ is a torus in $O(q)$, since this subgroup is isomorphic to the translation group of the torus.

Step 2. Let $G \subset O(q)$ be a toral subgroup such that $\text{Im}f=\text{Orb}_{G}(f(0))$. $\forall g \in G$, the map $f(x) \mapsto g \cdot f(x)$ is an isometry of $\text{Im}f$, meaning $G \subset \text{Isom}(\text{Im}f)$. The identity component of the Lie group $\text{Isom}(\text{Im}f)$ consists of parallel translations, i.e., $\text{Isom}_0(\text{Im}f) = \mathbb{T}^n$. By the definition of a toral group, $G$ is connected, hence $G \subset \mathbb{T}^n$. But then, since $\text{Im}f$ lies in the orbit of $G$, we get $\mathbb{T}^n=G$, which means the translation group of the torus extends to $O(q)$.

Step 3. From the previous step, we know that $\text{Im}f=\text{Orb}_{\mathbb{T}^n}(f(0))$, where $T^n \subset O(q)$ is a toral subgroup. Consider $\mathbb{T}_M$, a maximal torus in $O(q)$ containing $\mathbb{T}^n$. In some basis (in which, in particular, elements with indices $2j-1$ and $2j$ for $j=1,\dots,m$ must form an orthonormal basis of a plane), the elements of $\mathbb{T}_M$ have a special block-diagonal form. In the case $q=2m$, the form is as follows:

\[
\begin{pmatrix}
\cos{\xi_1} & -\sin{\xi_1} & \cdots & 0 & 0 \\
\sin{\xi_1} & \cos{\xi_1} & \cdots & 0 & 0\\
\vdots & \vdots & \ddots & \vdots & \vdots\\
0 & 0 & \cdots & \cos{\xi_m} & -\sin{\xi_m} \\
0 & 0 & \cdots & \sin{\xi_m} & \cos{\xi_m} 
\end{pmatrix}
\]

And if $q=2m+1$, then:

\[
\begin{pmatrix}
\cos{\xi_1} & -\sin{\xi_1} & \cdots & 0 & 0 & 0 \\
\sin{\xi_1} & \cos{\xi_1} & \cdots & 0 & 0 & 0\\
\vdots & \vdots & \ddots & \vdots & \vdots & \vdots\\
0 & 0 & \cdots & \cos{\xi_m} & -\sin{\xi_m} & 0 \\
0 & 0 & \cdots & \sin{\xi_m} & \cos{\xi_m} & 0\\
0 & 0 & \dots & 0 & 0 & 1
\end{pmatrix}
\]

Since $\mathbb{T}^n$ is a subgroup of $\mathbb{T}_M$, in this basis its elements will have the following block form:

\[
\begin{pmatrix}
\cos{f_1(x_1,\dots,x_n)} & -\sin{f_1(x_1,\dots,x_n)} &&&&\\
\sin{f_1(x_1,\dots,x_n)} & \cos{f_1(x_1,\dots,x_n)} & &&&\\
&& \ddots & &\\
&&&\cos{f_n(x_1,\dots,x_n)} & -\sin{f_n(x_1,\dots,x_n)} \\
&&&\sin{f_n(x_1,\dots,x_n)} & \cos{f_n(x_1,\dots,x_n)}\\
&&&&&0&\\
&&&&&&\ddots\\
&&&&&&&0
\end{pmatrix}
\]

Where all $f_j$ are smooth and additive (as it is an embedding of a Lie group), hence $\mathbb{R}$-linear, and also $\cos{f_j}$, $\sin{f_j}$ are $\Lambda \mathbb{Z}^n$-periodic functions, where $\Lambda$ is the lattice of the Gromov torus under consideration. This means $f_j(x)=\big< r_j,x \big>$, and furthermore $\big< r_j,x+a \big>=2\pi z_a+\big< r_j,x \big>$, where $a=\Lambda z \in \Lambda \mathbb{Z}^n$, $z_a \in \mathbb{Z}$. It follows that $\big< r_j,x \big>+\big< r_j,a \big>=2\pi z_a+\big< r_j,x \big>$, hence $\big< r_j,\Lambda z \big>=2\pi z_a$, which implies $\frac{1}{2\pi} \big< \Lambda^Tr_j,z \big> \in \mathbb{Z}$ for any $z \in \mathbb{Z}^n$. This yields $r_j=2\pi \Lambda^* k_j$ for some $k_j \in \mathbb{Z}^n$.

Let us look at how our map $f$ appears in these coordinates. If $x \in \mathbb{T}^n$, then $f(x)=x\cdot f(0)$, the matrix $x$ has the form described above, and $f(0) = (a_1,b_1, \dots,a_m,b_m)$. Then:

$$f(x)=(a_1\cos{2\pi \big< \Lambda^*k_1,x \big>}-b_1\sin{2\pi \big< \Lambda^*k_1,x \big>},a_1\sin{2\pi \big< \Lambda^*k_1,x \big>}+b_1\cos{2\pi \big< \Lambda^*k_1,x \big>},\dots)$$

Let $\{ \varepsilon_j \}_{j=1}^q$ be the chosen basis. Construct a basis $\{ e_j \}_{j=1}^q$ from it as follows: in each plane $\text{Span}(\{ \varepsilon_{2j-1}, \varepsilon_{2j}\})$ for $j=1,\dots,m$, the vectors $\varepsilon_{2j-1}, \varepsilon_{2j}$ differ from the vectors $e_{2j-1},e_{2j}$ by a rotation, and $f(0)$ has zero even coordinates in the new basis (in the case of an odd $q$, we do nothing to the last coordinate). Then, after renaming, the map $f$ in the new basis will take the following form:

$$f(x)=(a_1\cos{2\pi \big< \Lambda^*k_1,x \big>},a_1\sin{2\pi \big< \Lambda^*k_1,x \big>},\dots) = \sum_{j=1}^n{a_j\cos{(2\pi \big< \Lambda^*k_j,x \big>)}e_{2j-1}+a_j\sin{(2\pi \big< \Lambda^*k_j,x \big>)}e_{2j}}=$$

$$=\sum{\frac{a_j}{2}((u_j(x)+u_{-j}(x))e_{2j-1}-i(u_j(x)-u_{-j}(x))e_{2j})}=\sum{\frac{a_j}{2}(e_{2j-1}-ie_{2j})u_j(x)}+\frac{a_j}{2}(e_{2j-1}+ie_{2j})u_{-j}(x) =:$$

$$=: \sum{c_ju_j(x)+c_{-j}u_{-j}(x)}$$

Where $u_j(x)=\exp{(2\pi i \big< \Lambda^*k_j,x \big>)}$, $u_{-j}(x)=\exp{(-2\pi i \big< \Lambda^*k_j,x \big>)}$. Thus, we have brought $f$ into the form of its Fourier expansion. Here it is easy to see that $c_j \perp c_{-j}$ and $c_j \perp c_i$ for $i \neq j$.

Step 4. Let $0 \neq c_k = a_k+ib_k$, $a_k,b_k \in \mathbb{R}^q$. Since $\text{Im}f \subset \mathbb{R}^q$, we have $c_{-k}=\overline{c_k}$. Since $|k| = \frac{\sqrt{n}}{2\pi} \neq 0$, we have $k \neq -k$, and therefore $0=\big< c_k,c_{-k} \big> = \big< a_k+ib_k,a_k-ib_k \big> = |a_k|^2-|b_k|^2+2i\big< a_k,b_k \big>$. That is, $|a_k|=|b_k|$ and $a_k \perp b_k$.

Consider $c_l = a_l+ib_l$, $l \neq k$, $l \neq -k$ (if any exist). Then $c_l \perp c_k+c_{-k}$, $c_l \perp c_k-c_{-k}$, $c_{-l} \perp c_k+c_{-k}$, $c_{-l} \perp c_k-c_{-k}$, which implies $a_l \perp a_k$, $b_l \perp a_k$, $a_l \perp b_k$, $b_l \perp b_k$. Thus, $\{ \frac{a_k}{|a_k|}, \frac{b_k}{|b_k|} \} =: \{ e_j \}$ is an orthonormal system. Let $r_j := |a_j| = |b_j|$. Then we have:

$$f(x) = \sum_j{r_j(e_{2j-1}+ie_{2j})\text{exp}(2\pi i \big< k_j,x \big>)+r_j(e_{2j-1}-ie_{2j})\text{exp}(-2\pi i \big< k_j,x \big>)} = $$

$$=\sum{2r_j\cos{(-2\pi \big< k_j,x \big>)}e_{2j-1}+2r_j\sin{(-2\pi \big< k_j,x \big>)}e_{2j}} =: A(x)v$$

Where the vector $v=(2r_1,0,2r_2,0,\dots,2r_n,0,\dots,0)$, and $A(x)$ is the following block matrix:

\[
\begin{pmatrix}
\cos{(-2\pi \big< k_1,x \big>)} & -\sin{(-2\pi \big< k_1,x \big>)} &&&&\\
\sin{(-2\pi \big< k_1,x \big>)} & \cos{(-2\pi \big< k_1,x \big>)} & &&&\\
&& \ddots & &\\
&&&\cos{(-2\pi \big< k_n,x \big>)} & -\sin{(-2\pi \big< k_n,x \big>)} \\
&&&\sin{(-2\pi \big< k_n,x \big>)} & \cos{(-2\pi \big< k_n,x \big>)}\\
&&&&&0&\\
&&&&&&\ddots\\
&&&&&&&0
\end{pmatrix}
\]

That is, $\text{Im}f$ is the orbit of the toral subgroup of $O(q)$ consisting of matrices as above for $x \in \mathbb{T}^n$.

\begin{flushright}
$\blacksquare$
\end{flushright}

\quad

\textbf{Proposition 5.3.} There exists a non-equivariant embedded Gromov torus.

\quad

\textbf{Proof.} Let $\Lambda\subset\mathbb{R}^n$ be a lattice, and let $\Lambda^*$ be its dual lattice. Consider a finite set $\{ k_j \} \subset \Lambda^*$ such that $|k_j|=\frac{\sqrt{n}}{2\pi}$. Since $|k_j|=|-k_j|$, we denote $k_{-j}:=-k_j$. Then this set is indexed by $M:=\{ \pm1,\pm2,\dots,\pm m \}$. Consider the map:
\[
f\colon\mathbb{R}^n\to\mathbb{C}^{q+1}, \qquad
f(x)=\sum_{j\in M} c_j\,e^{2\pi i\langle k_j,x\rangle}.
\]
Where $\{ c_j \}$ is a set of variable complex vectors of variable dimension. Let us also set $c_{-j}=\overline{c_j}$, so that the map becomes real-valued.

As discussed earlier, for $f$ to define a Gromov torus, the following conditions are sufficient:
\begin{enumerate}
  \item $|f|\equiv 1$;
  \item $|D_{e_\alpha} f|\equiv 1$ for all $\alpha=1,\dots,n$;
  \item $\bigl\langle D_{e_\alpha} f,D_{e_\beta} f\bigr\rangle\equiv 0$ for $\alpha\neq\beta$;
  \item $|D^2_{e_\alpha} f|\equiv \sqrt{\dfrac{3n}{n+2}}$ for all $\alpha=1,\dots,n$;
  \item $|D^2_{\frac{e_\alpha+e_\beta}{\sqrt{2}}}f|=\sqrt{\frac{3n}{n+2}}$ for $\alpha \neq \beta$.
\end{enumerate}
Here $(e_1,\dots,e_n)$ is the standard basis in $\mathbb{R}^n$. We wish to construct injective (modulo lattice) $f$ under these conditions such that not all $\big< c_i,c_j \big>$ vanish for $i \neq j$.

Substituting the Fourier expansion $f(x)=\sum_{j\in M} c_j\,e^{2\pi i\langle k_j,x\rangle}$ into conditions (1)–(5), we obtain:

\begin{enumerate}

\item$$\sum_{j\in M}|c_j|^2
+\sum_{\substack{i,j\in M\\i\neq j}}
\bigl\langle c_{i},c_{j}\bigr\rangle
\,e^{2\pi i\langle k_i-k_j,x\rangle}\equiv1;$$

\item $$\sum_{j\in M}|c_j|^2(k_j^\alpha)^2
+\sum_{\substack{i,j\in M\\i\neq j}}
\bigl\langle c_{i},c_{j}\bigr\rangle
k_i^\alpha k_j^\alpha
\,e^{2\pi i\langle k_i-k_j,x\rangle}\equiv\frac{1}{4\pi^2};$$

\item $$\sum_{j\in M}|c_j|^2k_j^\alpha k_j^\beta
+\sum_{\substack{i,j\in M\\i\neq j}}
\bigl\langle c_{i},c_{j}\bigr\rangle
k_i^\alpha k_j^\beta
\,e^{2\pi i\langle k_i-k_j,x\rangle}
\equiv 0, \alpha \neq \beta;$$

\item $$\sum_{j\in M}|c_j|^2(k_j^\alpha)^4
+\sum_{\substack{i,j\in M\\i\neq j}}
\bigl\langle c_{i},c_{j}\bigr\rangle
(k_i^\alpha k_j^\alpha)^2
\,e^{2\pi i\langle k_i-k_j,x\rangle}\equiv\frac{3n}{16\pi^4(n+2)};$$

\item $$\sum_{j \in M}{|c_j|^2(k_j^{\alpha}+k_j^{\beta})^4}+\sum_{\substack{i,j \in M\\i \neq j}}{\big< c_i,c_j \big> (k_i^\alpha + k_i^\beta)^2(k_j^\alpha + k_j^\beta)^2\,e^{2\pi i\langle k_i-k_j,x\rangle}} \equiv \frac{3n}{4\pi^4(n+2)},\alpha \neq \beta.$$

\end{enumerate}

Where $k_j^\alpha$ is the $\alpha$-th coordinate of $k_j$ in the basis $\{e_\alpha\}$.

Notice that $\{ e^{2\pi i\langle k_i-k_j,x\rangle} \}$ for distinct differences $k_i-k_j$ form a part of the Fourier basis on the torus. Hence, functions (constants in this case) must expand in terms of them uniquely.

For convenience, let us denote $L := M\times M$, $L_v := \{(i,j)\in L\mid k_i-k_j=v\}, v\in\mathbb{R}^n$, and $2l_v:=|L_v|$ (each such class contains an even number of elements, since if $(i,j) \in L_v$, then $(-j,-i)\in L_v$). Then $L = \bigsqcup_{v \in \mathbb{R}^n}L_v$, the number of non-empty sets $L_v$ is finite, and $|L|=2m(2m-1)=\sum_{v \in \mathbb{R}^n}{2l_v}$. These classes correspond to different elements of the Fourier basis in the equations under consideration.

Let us write the conditions on the coefficients for $\{ e^{2\pi i\langle k_i-k_j,x\rangle} \}$ in the equations above. For $L_0$, we have:

\begin{enumerate}

\item$$\sum_{j\in M}|c_j|^2=1;$$

\item $$\sum_{j\in M}|c_j|^2(k_j^\alpha)^2
=\frac{1}{4\pi^2};$$

\item $$\sum_{j\in M}|c_j|^2k_j^\alpha k_j^\beta
=0, \alpha \neq \beta;$$

\item $$\sum_{j\in M}|c_j|^2(k_j^\alpha)^4
=\frac{3n}{16\pi^4(n+2)};$$

\item $$\sum_{j \in M}{|c_j|^2(k_j^{\alpha}+k_j^{\beta})^4} = \frac{3n}{4\pi^4(n+2)}.$$

\end{enumerate}

And for $L_v$, $v \neq 0$:

\begin{enumerate}

\item$$\sum_{(i,j)\in L_v}
\bigl\langle c_{i},c_{j}\bigr\rangle
=0;$$

\item $$\sum_{(i,j)\in L_v}
\bigl\langle c_{i},c_{j}\bigr\rangle
k_i^\alpha k_j^\beta
=0;$$

\item $$\sum_{(i,j)\in L_v}
\bigl\langle c_{i},c_{j}\bigr\rangle(k_i^\alpha + k_i^\beta)^2(k_j^\alpha + k_j^\beta)^2
=0.$$

\end{enumerate}

We can treat this system of equations as a linear system over $\mathbb{C}$ with respect to the variables $\{\big< c_i,c_j \big>\}$. For this, we must also add the following equations, otherwise the obtained solutions will not define valid complex inner products:

\begin{enumerate}

\item $$\big< c_i,c_j \big>=\overline{\big< c_j,c_i \big>}$$

\item $$\big< c_i,c_j \big> = \big< c_{-j},c_{-i} \big>$$

\end{enumerate}

Now, finally, we describe the construction of our non-equivariant Gromov torus. First, let us take an arbitrary Gromov torus $\mathbb{T}^n$ ($n \geq 3$), the existence of which was proven in \cite{MG}. Let $\Lambda$ be its corresponding lattice. Denote $\{k_j\}_{j=1}^m:=\frac{\sqrt{n}}{2\pi}\mathbb{S}^{n-1} \cap \Lambda^*$. Consider the system as above with coefficients depending on $\{ k_j \}$, corresponding to this torus. Notice that the entire system decouples into independent subsystems for the sets of variables $\{ \big< c_i,c_j \big>, (i,j) \in (L_v \cup L_{-v}) \}$. 

Let us first look at the subsystem corresponding to the variables $\{ |c_j|^2 \}$. Since the entire system must define a Gromov torus, this subsystem must have some non-negative solution, and all components of the solution cannot vanish simultaneously. Also, since $c_{-j} = \overline{c_j}$, we can restrict our view to variables with positive indices, which divides the right-hand side by two. Let us denote:

$$\Phi(k):=(1,(k^1)^2, \dots , (k^n)^2, k^1 \cdot k^2, \dots, k^{n-1}\cdot k^n, (k^1)^4, \dots, (k^n)^4,(k^1+k^2)^4,\dots,(k^{n-1}+k^n)^4)$$

$$b:=\bigg(\frac{1}{2}, \frac{1}{8\pi^2},\dots,\frac{1}{8\pi^2},0,\dots,0,\frac{3n}{32\pi^4(n+2)},\dots,\frac{3n}{32\pi^4(n+2)},\frac{3n}{8\pi^4(n+2)},\dots,\frac{3n}{8\pi^4(n+2)}\bigg)$$

That is, $\Phi:\frac{\sqrt{n}}{2\pi}\mathbb{S}^{n-1} \rightarrow \mathbb{R}^{2n^2+1}$ is a certain map, $\Phi(k_j)$ are the columns of the coefficient matrix of our subsystem, and $b$ is its right-hand side. Let $\{ \lambda_j \}_{j=1}^m$ be the solution of the subsystem corresponding to the Gromov torus, and let $\{ \lambda_j \}_{j=1}^r$ be the positive components of this solution, so that $\sum_{j=1}^r{\lambda_j\Phi(k_j)}=b$. Denote $V_0:=\{ \Phi(k_1), \dots, \Phi(k_r) \}$.

Let $\varphi \in O(n)$ be an element of the isometry group of $\frac{\sqrt{n}}{2\pi}\mathbb{S}^{n-1}$. Note that the lattice $\varphi(\Lambda)$ corresponds to a Gromov torus isometric to the original one, with $\Phi(\varphi(k_j))$ being the columns of the corresponding subsystem, while the right-hand side remains equal to $b$. Then this system will have a non-negative solution. In other words, there exist $\mu_j \geq 0$ such that $\sum_{j=1}^l{\mu_j\Phi(\varphi(k_j))}=b$.

First, assume that $\text{Span}(V_0)$ does not coincide with $\text{Span}(\Phi(\frac{\sqrt{n}}{2\pi}\mathbb{S}^{n-1}))$. We prove that in this case, there exists $\varphi_1 \in O(n)$ such that $\Phi(\varphi_1(k_j)) \notin \text{Span}(V_0)$ for $j=1,\dots,m$. Note that $\Phi$ is an analytic map, as it is the restriction of a polynomial map from $\mathbb{R}^n$ to the analytic submanifold $\frac{\sqrt{n}}{2\pi}\mathbb{S}^{n-1}$. In particular, it is continuous. $\text{Span}(V_0)$ is a closed set, so $K:=\Phi^{-1}(\text{Span}(V_0))$ is also closed in the sphere, meaning it is compact. Also notice that $\text{Int}(K) = \emptyset$ in the topology of the sphere, because otherwise $\Phi(\frac{\sqrt{n}}{2\pi}\mathbb{S}^{n-1}) \subset \text{Span}(V_0)$ by the uniqueness theorem for analytic functions, which contradicts our assumption. Hence $\lambda_{\mathbb{S}}^{n-1}(K)=0$, where $\lambda_{\mathbb{S}}^{n-1}$ is the Lebesgue measure on $\frac{\sqrt{n}}{2\pi}\mathbb{S}^{n-1}$.

Consider $O(n)$ as a Lie group with a left-invariant Haar measure $\mu$. Denote $G_j:O(n) \rightarrow \frac{\sqrt{n}}{2\pi}\mathbb{S}^{n-1}$ by $G_j(g):=g(k_j)$, where $j=1,\dots,m$. Note that $G_j$ is a smooth and surjective map for all $j$. Also, each $G_j$ induces a certain Borel measure on the sphere. Namely, if $A$ is a Borel set in the topology of $\frac{\sqrt{n}}{2\pi}\mathbb{S}^{n-1}$, let us denote $\lambda(A):=\mu(G_j^{-1}(A))$. Notice that if $h \in O(n)$, then $\lambda(hA)=\mu(G_j^{-1}(hA))=\mu(hG_j^{-1}(A))=\mu(G_j^{-1}(A))=\lambda(A)$. Thus, $\lambda$ is a Borel measure invariant under $O(n)$, and therefore it is proportional to the Lebesgue measure $\lambda_{\mathbb{S}}^{n-1}$. From this, we obtain $0=\lambda(K)=\mu(\{ g \in O(n), g(k_j) \in K \})$. Consequently, $\mu(g\in O(n),\exists j=1,\dots,l:g(k_j)\in K)=0$. Thus, almost any element of $O(n)$ can be chosen as $\varphi_1$.

Now denote $m_j^1:=\Phi(\varphi_1(k_j))$. By the arguments above, there exist $\mu_j^1 \geq 0$ such that $\sum_{j=1}^l{\mu_j^1m_j^1}=b$. Without loss of generality, let $\mu_1^1, \dots, \mu_{r_1}^1 > 0$ (reindexing if necessary). Denote $V_1:=V_0 \cup\{m_1^1,\dots,m_{r_1}^1\}$. Note that by construction $\text{dim}(\text{Span}(V_1)) > \text{dim}(\text{Span}(V_0))$.

If $\text{dim}(\text{Span}(V_1))$ is still not maximal, we can similarly construct $\varphi_2 \in O(n)$, $m_j^2:=\Phi(\varphi_2(k_j))$, $j=1,\dots,l$, such that $m_j^2 \notin \text{Span}(V_1)$. In this case, there exist $\mu_1^2, \dots, \mu_{r_2}^2>0$ such that $\sum_{j=1}^{r_2}{\mu_j^2m_j^2}=b$. Denoting $V_2:=V_1 \cup \{m_1^2,\dots,m_{r_2}^2\}$, we have by construction that $\text{dim}(\text{Span}(V_2)) > \text{dim}(\text{Span}(V_1))$.

Proceeding in this manner, sooner or later, because the dimension of $V_j$ increases each time, we can ensure that at step $R$ we have $V_R = \text{Span}(\Phi(\frac{\sqrt{n}}{2\pi}\mathbb{S}^{n-1}))$. That is, $V_R$ is already a spanning set for this vector subspace. Moreover, by construction:

$$\sum_{j=1}^r\frac{\lambda_j}{R}\Phi(k_j)+\sum_{j=1}^{r_1}\frac{\mu^1_j}{R}m_j^1+\dots+\sum_{j=1}^{r_R}\frac{\mu^R_j}{R}m_j^R=b$$

Where all $\frac{\lambda_j}{R},\frac{\mu_j^i}{R}>0$. For uniformity, let us rename all vectors and coefficients in the equation above: $\sum_{j=1}^{\rho}\lambda_jm_j=b$, where all $\lambda_j>0$, and $\{ m_j \}_{j=1}^{\rho}$ form $\text{Span}(\Phi(\frac{\sqrt{n}}{2\pi}\mathbb{S}^{n-1}))$, and let $\Phi(k_j)=m_j$, which exist by construction.

Now, using the construction above, we will modify the lattice $\Lambda$. Consider points $\tilde{k_j} \in \Phi(\frac{\sqrt{n}}{2\pi}\mathbb{S}^{n-1})$ such that $\tilde{k_j}$ is close to $k_j$ for each $j$, and $\frac{\tilde{k_j}}{\frac{\sqrt{n}}{2\pi}}$ are rational points in $\mathbb{R}^n$ (hereinafter, for convenience, we will call this property the rationality of points on such a sphere). Since $\text{Span}(m_1,\dots,m_{\rho})=\text{Span}(\Phi(\frac{\sqrt{n}}{2\pi}\mathbb{S}^{n-1}))$ and $\sum_{j=1}^{\rho}\lambda_jm_j=b$ for some $\lambda_j > 0$, $b$ lies strictly inside the cone $\{ a_1m_1+\dots+a_{\rho} m_{\rho} \mid a_1, \dots a_{\rho} \geq 0\} \subset \text{Span}(\Phi(\frac{\sqrt{n}}{2\pi}\mathbb{S}^{n-1}))$, and this cone has a non-empty interior relative to $\text{Span}(\Phi(\frac{\sqrt{n}}{2\pi}\mathbb{S}^{n-1}))$. Then, if $\tilde{k_j}$ are sufficiently close to $k_j$, and consequently $\Phi(\tilde{k_j})$ are sufficiently close to $m_j$, $b$ will still lie strictly inside the cone $\{ a_1\Phi(\tilde{k_1})+\dots+a_{\rho} \Phi(\tilde{k_{\rho}}) \mid a_1, \dots a_{\rho} \geq 0\}$ with a non-empty interior relative to $\text{Span}(\Phi(\frac{\sqrt{n}}{2\pi}\mathbb{S}^{n-1}))$.

Now, via the transformation above, we can assume that all $k_j$ lie on some sufficiently small orthogonal lattice $\Lambda^*$. For now, we do not consider vectors on this lattice other than our $k_j$, with respect to which the system has a strictly positive solution. It is claimed that now, when adding any column $\Phi(k)$ for $k \in \frac{\sqrt{n}}{2\pi}\mathbb{S}^{n-1}$ to our system, it will still have a strictly positive solution. Indeed, for the solution to be positive, it suffices to find $\tilde{\lambda_j}>0$ and $t>0$ such that $\sum_{j=1}^{\rho}\tilde{\lambda_j}m_j=b-t\Phi(k)$. Since $b$ lies strictly within the interior of the cone under consideration, for a sufficiently small positive $t$, the point $b-t\Phi(k)$ will also lie in the interior, meaning that the corresponding positive $\tilde{\lambda_j}$ will exist. Thus, we can assume that $\{ k_j\}$ are all the points from the set $\frac{\sqrt{n}}{2\pi}\mathbb{S}^{n-1} \cap \Lambda^*$. Moreover, by making the lattice sufficiently small, we can add any finite set of rational (up to scaling) vectors $k_j$ into it. We will exploit this idea below.

Now let us look at the subsystems corresponding to the variables $\{ \big< c_i,c_j \big>, (i,j) \in (L_v \cup L_{-v}) \}$ for $v \neq 0$. Obviously, satisfying the subsystem for $L_v$ implies satisfying it for $L_{-v}$ as well, since the values of the variables in them differ only by conjugation. Therefore, we may speak only of the subsystem for $L_v$, and without the condition $\big< c_i,c_j \big> = \big< c_{-j},c_{-i} \big>$. If in some subsystem corresponding to a non-empty class $L_v$ ($v \neq 0$) there are more variables than equations, it will possess a non-trivial solution since it is homogeneous. This would allow us to construct a set of inner products satisfying all the conditions above, such that not all $\big< c_i,c_j \big>$ vanish for $i \neq j$. In the system for a non-empty class $L_v$, we have $2n^2+1+l_v$ equations and $2l_v$ unknowns.

We show that if we take a rectangular lattice that is sufficiently small, the condition above will be met. Consider $\frac{\sqrt{n}}{2\pi}\mathbb{S}^2 \subset \{ x \in\mathbb{R}^n, x_4 = \dots = x_n = 0 \}$. Let $S$ be some non-major circle of this sphere lying in a plane perpendicular to the $Ox_1$ axis and intersecting this axis at a point with a coordinate of the form $(\frac{\sqrt{n}}{2\pi}q,0,0)$, where $q$ is some positive rational number. Then there will be infinitely many rational (in our sense) points on this circle. Consider $2(2n^2+2)$ non-intersecting pairs of rational points on it, symmetric with respect to $Ox_1$. It is easy to see that if $(u,v)$ is one of such pairs, then $u+v=(2\frac{\sqrt{n}}{2\pi}q,0,0)$. Now, using the previously described algorithm, we can assume that for any such pair: $u,-v \in \Lambda^*$, and then $|L_{(2\frac{\sqrt{n}}{2\pi}q,0,0)}| \geq 2(2n^2+2)$, which means $l_{(2\frac{\sqrt{n}}{2\pi}q,0,0)} > 2n^2+1$.

Let us redefine our lattice one final time. Suppose that $\{ k_j \}$ is the set of all rational points on the sphere we have chosen so far. By construction, we may assume that there exist $n$ $\mathbb{R}$-lineary independent points among them. Define $\Gamma^*:=\{ \sum{z_j \cdot k_j} \big| z_j \in \mathbb{Z} \}$. This is an infinite discrete subgroup of $\Lambda^*$, and every point in $\mathbb{R}^n$ is within a bounded distance from $\Gamma^*$. Hence, it is a lattice of rank $n$. We then define $\Gamma := (\Gamma^*)^*=\{ v \in \mathbb{R}^n \big| \forall k \in \Gamma^*: \big< k,v \big> \in \mathbb{Z} \}$. This will be our sought-after lattice.

Thus, we have constructed a lattice $\Gamma$ and a set of complex numbers $\{\big< c_i, c_j \big>\}_{i,j\in M}$ that satisfies the corresponding system. Here, all $|c_j|^2 > 0$, and $\exists\,i\neq j:\ \langle c_i,c_j\rangle\neq0$. Now we want to reconstruct the vectors $\{ c_j \}$ themselves, such that obtained $f$ is injective modulo $\Gamma$. It would then follow from Lemmas 5.1 and 5.2 that we have constructed a non-equivariant embedded Gromov torus for the lattice $\Gamma$.

Note that since the system of equations for the numbers $\{\langle c_i,c_j\rangle\}_{i,j\in M,\;i\neq j}$ is homogeneous, we can multiply them all simultaneously by the same constant. In particular, we can make them arbitrarily small in absolute value.

It is convenient to write the inner products as follows. Let $c_i=a_i+ib_i$, where $a_i,b_i\in\mathbb{R}^q$. Then:
\begin{align*}
\langle c_i,c_j\rangle
 &= \langle a_i,a_j\rangle-\langle b_i,b_j\rangle
   + i\bigl(-\langle a_i,b_j\rangle+\langle a_j,b_i\rangle\bigr),\\
\langle c_i,c_{-j}\rangle
 &= \langle a_i,a_j\rangle+\langle b_i,b_j\rangle
   + i\bigl(\langle a_i,b_j\rangle+\langle a_j,b_i\rangle\bigr).
\end{align*}
It suffices to find the real vectors $\{a_i,b_i\}$; from them, we can reconstruct $c_i$. Also, from the known numbers $\{\langle c_i,c_j\rangle\}$, we know all the quantities $\{\langle a_i,b_j\rangle,\ \langle a_i,a_j\rangle,\ \langle b_i,b_j\rangle\}$ by virtue of the identities given above.

For convenience, we combine all $a_i$ and $b_i$ into a single system and denote them by $\{u_s^\varepsilon\}_{s=1}^N, u_s^\varepsilon\in\mathbb{R}^n$, without distinguishing between $a_i$ and $b_i$, where $\varepsilon$ is the maximum modulus of the varied inner products. We are given the lengths $|u_j^\varepsilon|$ and the inner products $\langle u_i^\varepsilon,u_j^\varepsilon\rangle =: s_{ij}(\varepsilon)$, where for $i\neq j$ we have $s_{ij}(\varepsilon)\to0$ as $\varepsilon\to0$, since we can make all $\langle c_i,c_j\rangle$ (for $i\neq j$) arbitrarily small in modulus. Notice that in this case $\angle(u_i^\varepsilon,u_j^\varepsilon)\to\frac{\pi}{2}$ as $\varepsilon\to0$.

We prove by induction that for any sufficiently small $\varepsilon$, there exists a set of linearly independent vectors $\{ u_j^{\varepsilon} \}_{j=1}^n \subset \mathbb{R}^n$ that realize our constructed conditions on their non-zero lengths and inner products depending on the parameter $\varepsilon$. 

\emph{Base case: $n=2$.} Taking $\varepsilon$ small enough to satisfy the strict Cauchy-Bunyakovsky-Schwarz inequality, we can construct a pair of vectors satisfying the given conditions.

\emph{Inductive step: $n\mapsto n+1$.} Suppose that for some $n$ we have already constructed $u_1^\varepsilon,\dots,u_n^\varepsilon\in\mathbb{R}^n$ with the required conditions. We want to add a vector $u_{n+1}^\varepsilon$ with a fixed non-zero length, for which:
\[
\begin{cases}
\langle x,u_1^\varepsilon\rangle = d_1(\varepsilon),\\
\langle x,u_2^\varepsilon\rangle = d_2(\varepsilon),\\
\quad\vdots\\
\langle x,u_n^\varepsilon\rangle = d_n(\varepsilon),
\end{cases}
\]
where $d_j(\varepsilon)\to0$ as $\varepsilon\to0$ (these are determined from the conditions on the inner products with $u_{n+1}^\varepsilon$). Our system can be written as $A(\varepsilon)x = d(\varepsilon)$, where the columns of the matrix $A(\varepsilon)$ are the vectors $u_j^{\varepsilon}$. Since $\angle(u_i^{\varepsilon}, u_j^{\varepsilon}) \to \frac{\pi}{2}$ as $\varepsilon \to 0$ and the lengths $|u_j^{\varepsilon}|$ are constant and non-zero, the matrix $A(\varepsilon)$ converges to a matrix whose columns are orthogonal non-zero vectors. In particular, there exists an $\varepsilon_0 > 0$ such that $A(\varepsilon)$ is non-singular for all $\varepsilon < \varepsilon_0$. By Cramer's rule, matrix inversion is continuous on $\text{GL}_n(\mathbb{R})$, so $||A(\varepsilon)^{-1}||$ is bounded as $\varepsilon \to 0$.

Hence, for sufficiently small $\varepsilon$, we have $x = A(\varepsilon)^{-1}d(\varepsilon)$. Since $|d(\varepsilon)| \to 0$ and $||A(\varepsilon)^{-1}||$ is bounded, $|x| \to 0$ as $\varepsilon \to 0$.

Let us add one more coordinate to our space, orthogonal to $\mathbb{R}^n$. Now we have $\mathbb{R}^{n+1}$, and we want to find $u_{n+1}=x+y$, with $|u_{n+1}|^2=|x|^2+|y|^2$ for some $y \perp \mathbb{R}^n$, where $|u_{n+1}|$ is fixed. Since $|x| \rightarrow 0$, such a $y$ exists for sufficiently small $\varepsilon$. Thus, for any sufficiently small $\varepsilon$, we have constructed $u_{n+1}$ satisfying all the required conditions.

Finally, we observe that the map $f$ is injective modulo $\Gamma$. Suppose that $f(x) = f(y)$. Since all $c_j$ are linearly independent, it follows that $\langle k_j, x - y \rangle \in \mathbb{Z}$ for all $j$. Because $\Gamma^*$ is generated by $\{ k_j \}$ over $\mathbb{Z}$, we have $\langle k, x - y \rangle \in \mathbb{Z}$ for all $k \in \Gamma^*$. This implies $x - y \in \Gamma$. Thus, $f$ is the required map.

\begin{flushright}
$\blacksquare$
\end{flushright}

\quad

\textbf{Remark}. Notice that, according to our construction and the reconstruction of $\{ c_j \}$, any such torus has as many continuous parameters for its embedding as there are homogeneous subsystems with non-trivial solutions in its system. Since we can create as many homogeneous subsystems with non-trivial solutions as desired, we can actually construct a non-equivariant embedded Gromov torus whose minimal embeddings have arbitrarily many continuous parameters.

\quad

\textbf{Acknowledgments.} I would like to thank my academic supervisor, Nina Lebedeva, for her continuous guidance. I also wish to thank Anton Petrunin for suggesting the initial problem and overall help, and Robert Bryant for his valuable comments on this topic. This paper was prepared as a course project at SPbU's MCS faculty.



\end{document}